# DISORDERED PINNING MODELS AND COPOLYMERS: BEYOND ANNEALED BOUNDS

By Fabio Lucio Toninelli [1]

*École Normale Supérieure de Lyon and CNRS*

We consider a general model of a disordered copolymer with adsorption. This includes, as particular cases, a generalization of the copolymer at a selective interface introduced by Garel et al. [*Europhys. Lett.* **8** (1989) 9–13], pinning and wetting models in various dimensions, and the Poland–Scheraga model of DNA denaturation. We prove a new variational upper bound for the free energy via an estimation of noninteger moments of the partition function. As an application, we show that for strong disorder the quenched critical point differs from the annealed one, for example, if the disorder distribution is Gaussian. In particular, for pinning models with loop exponent $0 < \alpha < 1/2$ this implies the existence of a transition from weak to strong disorder. For the copolymer model, under a (restrictive) condition on the law of the underlying renewal, we show that the critical point coincides with the one predicted via renormalization group arguments in the theoretical physics literature. A stronger result holds for a "reduced wetting model" introduced by Bodineau and Giacomin [*J. Statist. Phys.* **117** (2004) 801–818]: without restrictions on the law of the underlying renewal, the critical point coincides with the corresponding renormalization group prediction.

**1. Introduction.** We consider a rather general class of directed polymers interacting with a one-dimensional defect through a quenched disordered potential. This includes various models motivated by (bio-)physics: among others, wetting models in $(1 + 1)$ dimensions [12, 14], pinning of $(1+d)$-dimensional directed polymers on columnar defects [27], copolymers at selective interfaces [16, 25] and the Poland–Scheraga (PS) model of DNA denaturation [11, 24]. For further motivations and references, we refer to [17], Chapter 1. One of the interesting aspects of these models is that they

Received September 2007; revised September 2007.
[1]Supported in part by the GIP-ANR project JC05_42461 (POLINTBIO).

present a nontrivial localization–delocalization phase transition due to an energy–entropy competition.

Mathematically, the model is defined in terms of a renewal sequence whose inter-arrival law has a power-like tail with exponent $\alpha + 1 \geq 1$. The model is exactly solvable in absence of disorder, and it turns out that the transition can be of any given order, from first to infinite, according to the value of $\alpha$. This is therefore an ideal testing ground for physical arguments (Harris criterion, renormalization group computations) and predictions concerning the effect of disorder on the critical exponents and on the location of the critical point.

The comprehension of this model has witnessed remarkable progress on the mathematical side, as proved by the recent book [17]. In particular it has been shown that for wetting, pinning or PS models, an arbitrary amount of disorder modifies the free-energy critical exponent if $\alpha > 1/2$ [18], that is, *disorder is relevant* in this case, in agreement with the predictions of the so-called Harris criterion [23]. On the other hand, for $0 < \alpha < 1/2$ it has been proven recently [2, 29] that if disorder is weak enough the free-energy critical exponent coincides with that of the homogeneous (i.e., nondisordered) model, and the (quenched) critical point coincides with the annealed one: *disorder is irrelevant* (again, in agreement with the Harris criterion). These results about "irrelevance" of disorder for $0 < \alpha < 1/2$ have been later refined and complemented in [21] with results about correlation-length critical exponents. The *marginal case* $\alpha = 1/2$ is strongly debated in the theoretical physics literature: Ref. [14] claims that quenched and annealed critical points coincide for disorder weak enough, while [12] concludes the opposite and gives a precise prediction for their difference. See [2] and [29] for rigorous results in the marginal case, which however do not solve the controversy.

Here we attack two major open problems:

- Do quenched and annealed critical points coincide for *strong disorder*? The Harris criterion, which is based on the analysis of the stability of the homogeneous model to the addition of weak randomness, makes no prediction about this point for pinning models with $\alpha \leq 1/2$ or for the copolymer with any $\alpha$. [Here and in the following, we say for brevity "pinning models" but we actually include wetting and PS models, besides pinning of $(1 + d)$-dimensional directed polymers on columnar defects. Mathematically all these are variants of the same model, cf. beginning of Section 3.1.]
- For the copolymer model, it is known that the critical point is bounded above by the annealed one and below by an $\alpha$-dependent expression found by non-rigorous renormalization group arguments [25]. Are either of these two bounds optimal?

In this work we prove a new upper bound on the free energy of the model which in some cases is sufficient to answer these two questions. In particular, consequences of our bound include the following:

1. Both for pinning and copolymer models with, say, Gaussian randomness, for large disorder quenched and annealed critical points differ. We would like to emphasize that, especially for the "marginal case" of the wetting model with loop exponent $\alpha = 1/2$, this question was subject to dispute even very recently [12, 15].
2. We identify the strong-disorder behavior of the critical point both for Gaussian pinning models and for a "reduced" wetting model introduced by Bodineau and Giacomin [5] (for the "reduced model," the same result was proven recently by Bolthausen, Caravenna and de Tilière [6]; their proof is however based on a very different method).
3. For the copolymer model we prove that, as soon as a homogeneous depinning term is present in the Hamiltonian, quenched and annealed critical points differ *for every strength of the disorder*; in particular, the much-studied *critical slope at the origin* (cf. Section 3.3) is in this case strictly smaller than 1.
4. Finally, again for the copolymer model we prove that, if the law of the underlying renewal sequence satisfies a certain explicit condition [cf. equation (3.43) below; in particular, the condition requires the renewal to be transient], the critical point predicted by nonrigorous renormalization group arguments is indeed the correct one. This is however believed *not* to be the case in general, that is, if (3.43) does not hold.

Our basic idea is to estimate noninteger moments $\mathbb{E}Z^\gamma$ of the partition function with $1/(1+\alpha) \leq \gamma < 1$. The reason why we cannot go down to $\gamma < 1/(1+\alpha)$ is *not just technical* and will become clear soon. The method of fractional moments has been already applied successfully to the study of other quenched disordered models: we would like to mention in particular (a) the work [1] by Aizenman and Molchanov, where bounds on small moments (of order less than 1) of the resolvent kernel of random Schrödinger operators are employed to prove the occurrence of Anderson localization for strong disorder or extreme energies; (b) Ref. [8] by Buffet, Patrick and Pulé who compute exactly the free energy of a directed polymer in random environment on a regular tree via an estimation of $\mathbb{E}Z^\gamma$ with $1 < \gamma < 2$; and (c) Ref. [13] where Derrida and Evans, again via an estimation of $\mathbb{E}Z^\gamma$ with $1 < \gamma < 2$, improve previously known estimates on the critical temperature of directed polymers in random environment on finite-dimensional lattices.

We would like to conclude this introduction with two remarks. First of all, for pinning models with $0 < \alpha < 1/2$ and Gaussian randomness our results, together with those of [2] or [29], imply that there is a nontrivial transition from a weak- to a strong-disorder regime, see Remark 3.3 below for a precise statement. Secondly, Theorem 3.6 together with the numerical simulations of Ref. [9] strongly indicates that the critical point of the copolymer depends not only on $\alpha$ but also on the details of the inter-arrival law of the renewal, a possibly nonintuitive fact.

**2. The model and the main result.** Let $\tau := \{\tau_0, \tau_1, \ldots\}$ be a renewal sequence on the integers, with inter-arrival law $K(\cdot)$: $\tau_0 = 0$ and $\{\tau_i - \tau_{i-1}\}_{i \geq 1}$ is a sequence of IID random variables taking values in $\mathbb{N} \cup \{+\infty\}$, with law $\mathbf{P}(\tau_1 = n) =: K(n)$. We assume that for $n \in \mathbb{N}$

$$K(n) = \frac{L(n)}{n^{1+\alpha}} \tag{2.1}$$

with $0 < \alpha < \infty$ and $L(\cdot)$ a function varying slowly at infinity, that is, a positive function such that $\lim_{x \to \infty} L(rx)/L(x) = 1$ for every $r > 0$. In general $\sum_{n \in \mathbb{N}} K(n) = 1 - \mathbf{P}(\tau_1 = +\infty) \leq 1$.

A very popular example in the literature is the one where $K(\cdot)$ is the law of the first return to zero of the one-dimensional simple random walk $\{S_k\}_{k \geq 0}$, that is, $K(n) = K^{SRW}(n) := \mathbf{P}(\inf\{k > 0 : S_k = 0\} = 2n | S_0 = 0)$ (in which case $\alpha = 1/2$, $\sum_{n \geq 1} K(n) = 1$ and $L(\cdot) \sim const$).

In order to make a quick connection with the (bio)-physical systems presented in the introduction, let us mention that in the case of $(1+d)$-dimensional pinning models one takes $\alpha = d/2 - 1$ if $d \geq 2$ and $\alpha = 1/2$ if $d = 1$, in $(1+1)$-dimensional wetting models and copolymers at a selective interface usually $\alpha = 1/2$, while $\alpha \simeq 1.15$ in the case of the PS model in three dimensions [24]. Actually, for the copolymer model and the pinning model with $d = 1$ one usually makes the specific choice $K(\cdot) = K^{SRW}(\cdot)$.

We consider random copolymers with adsorption. The system has size $N \in \mathbb{N}$, and it is characterized by the parameters $\beta \geq 0$, $\lambda \geq 0$, $h \in \mathbb{R}$ and $\widetilde{h} \geq 0$. Moreover, we let $\omega := \{\omega_n\}_{n \in \mathbb{N}}$ and $\widetilde{\omega} := \{\widetilde{\omega}_n\}_{n \in \mathbb{N}}$ be the realizations of two IID sequences of random variables. We assume $\omega$ to be independent of $\widetilde{\omega}$. The joint law of $(\omega, \widetilde{\omega})$ is denoted by $\mathbb{P}$.

The partition function of the model is

$$Z_{N,\omega,\widetilde{\omega}} := \mathbf{E}\left[ e^{\sum_{n=1}^{N}(\beta\omega_n + h)\mathbf{1}_{\{n \in \tau\}}} \right. \tag{2.2}$$
$$\left. \times \prod_{j=1}^{I_N(\tau)} \left( \frac{1 + e^{-2\lambda \sum_{n=\tau_{j-1}+1}^{\tau_j}(\widetilde{\omega}_n + \widetilde{h})}}{2} \right) \mathbf{1}_{\{N \in \tau\}} \right],$$

where $I_N(\tau) := |\tau \cap \{1, \ldots, N\}| = \sum_{n=1}^{N} \mathbf{1}_{\{n \in \tau\}}$, while the free energy is

$$F(\beta, h, \lambda, \widetilde{h}) := \lim_{N \to \infty} \frac{1}{N} \log Z_{N,\omega,\widetilde{\omega}}. \tag{2.3}$$

We assume as usual the existence of all exponential moments of $\omega_1$:

$$(2.4) \qquad M(u) := \mathbb{E}(e^{u\omega_1}) < \infty$$

for every $u \in \mathbb{R}$, and similarly for $\widetilde{M}(u) := \mathbb{E}(\exp(u\widetilde{\omega}_1))$, and we set by convention $\mathbb{E}\omega_1^2 = \mathbb{E}\widetilde{\omega}_1^2 = 1$ and $\mathbb{E}\widetilde{\omega}_1 = 0$ (this can always be achieved via a redefinition of $\beta$, $\lambda$ and $\widetilde{h}$). For later convenience (cf. Section 3.2), we *do not* assume in general that the $\omega_n$'s are centered, although this can always be obtained via a trivial shift of $h$. Under these assumptions on $\mathbb{P}$, it is well known that the limit in (2.3) exists almost surely and in $L^1(\mathbb{P})$, and that it is almost surely independent of $(\omega, \widetilde{\omega})$. Another well-established fact is that $F(\beta, h, \lambda, \widetilde{h}) \geq 0$, which is an immediate consequence of

$$(2.5) \qquad Z_{N,\omega,\widetilde{\omega}} \geq \frac{e^{\beta\omega_N + h}}{2} K(N)$$

and (2.1). For a proof and a discussion of these facts, see for instance [17], Chapters 1 and 4.

It is actually one of the main questions in this context to decide when the free energy vanishes and when it is positive. The reason is the following: One usually defines the *localized region* as $\mathcal{L} := \{(\beta, h, \lambda, \widetilde{h}) : F(\beta, h, \lambda, \widetilde{h}) > 0\}$ and the *delocalized region* as $\mathcal{D} := \{(\beta, h, \lambda, \widetilde{h}) : F(\beta, h, \lambda, \widetilde{h}) = 0\}$. Then it is well known that, if $(\beta, h, \lambda, \widetilde{h}) \in \mathcal{L}$, the ratio $\mathbf{E}_{N,\omega,\widetilde{\omega}}(I_N(\tau))/N$ converges for $N \to \infty$ to a positive limit, almost-surely independent of $(\omega, \widetilde{\omega})$ (where $\mathbf{E}_{N,\omega,\widetilde{\omega}}$ denotes the disorder-dependent Gibbs average), while the same limit is zero in the interior of $\mathcal{D}$. This explains the names given to the two regions. We refer to [19] and [17], Chapter 8, for more precise and more refined estimates on $I_N(\tau)$ in the interior of $\mathcal{D}$, to [4], [20] and [17], Chapter 7, for path properties of $\tau$ in $\mathcal{L}$ and finally to [28] for estimates on $I_N(\tau)$ at the boundary between $\mathcal{D}$ and $\mathcal{L}$.

A very cheap way of showing that the system is delocalized for given values of $(\beta, h, \lambda, \widetilde{h})$ is via the annealed upper bound on the free energy:

$$(2.6) \quad \begin{aligned} F(\beta, h, \lambda, \widetilde{h}) &\leq F^{ann}(\beta, h, \lambda, \widetilde{h}) \\ &:= \lim_{N\to\infty} \frac{1}{N} \log \mathbf{E}\Bigg[e^{I_N(\tau)(h + \log M(\beta))} \\ &\qquad\qquad \times \prod_{j=1}^{I_N(\tau)} \bigg(\frac{1 + e^{(-2\widetilde{h}\lambda + \log \widetilde{M}(-2\lambda))(\tau_j - \tau_{j-1})}}{2}\bigg) \mathbf{1}_{\{N \in \tau\}}\Bigg], \end{aligned}$$

which is simply Jensen's inequality: $\mathbb{E} \log Z \leq \log \mathbb{E} Z$. If we define the delocalized region of the annealed model as $\mathcal{D}^{ann} := \{(\beta, h, \lambda, \widetilde{h}) : F^{ann}(\beta, h, \lambda, \widetilde{h}) = 0\}$, one has then obviously $\mathcal{D}^{ann} \subseteq \mathcal{D}$.

The point of this work is to provide an improved upper bound on $F$ which is enough to conclude that $\mathcal{D}^{ann} \neq \mathcal{D}$. In some cases, we will even be able to identify sharply the boundary between the two regions.

Going beyond annealing has appeared so far to be a difficult task. A natural idea is to try *Morita-type bounds* [26], that is, constrained annealing. In other words, for every function $p(\omega, \widetilde{\omega})$ such that $\mathbb{E}p(\omega, \widetilde{\omega}) = 0$, it is easily seen that

$$(2.7) \qquad F(\beta, h, \lambda, \widetilde{h}) \leq \lim_{N \to \infty} \frac{1}{N} \log \mathbb{E}(e^{p(\omega, \widetilde{\omega})} Z_{N, \omega, \widetilde{\omega}}).$$

This can indeed improve the upper bound (2.6) on $F$ if $p$ is suitably chosen but, as shown in [10], it cannot improve the estimate $\mathcal{D}^{ann} \subseteq \mathcal{D}$, as long as $p$ is a sum of local functions ("finite-order Morita approximation"):

$$p(\omega, \widetilde{\omega}) = \sum_{n=1}^{N} p_0(\theta_n(\omega, \widetilde{\omega})),$$

with $p_0(\omega, \widetilde{\omega})$ a bounded function depending only on a finite number of $(\omega_i, \widetilde{\omega}_i)$, and $\theta$ being the shift operator $(\theta_n(\omega, \widetilde{\omega}))_m = (\omega_{n+m}, \widetilde{\omega}_{n+m})$. We will show that the estimation of noninteger moments of the partition function allows to bypass this difficulty.

In order to formulate our results we need a few auxiliary quantities. For $0 < \gamma \leq 1$ let

$$(2.8) \qquad c(\gamma) := \sum_{n \in \mathbb{N}} (K(n))^{\gamma} \geq c(1) = 1 - \mathbf{P}(\tau_1 = +\infty),$$

with strict inequality if $\gamma < 1$. We remark that $c(\cdot)$ is decreasing and that $c(\gamma) < \infty$ if $1/(1+\alpha) < \gamma \leq 1$. Also, observe that

$$c(1/(1+\alpha)) = \sum_{n \in \mathbb{N}} \frac{L(n)^{1/(1+\alpha)}}{n}$$

can be finite or infinite according to the behavior at infinity of $L(\cdot)$.

Define also, for $n \in \mathbb{N}$ and $\gamma$ such that $c(\gamma) < \infty$,

$$(2.9) \qquad \hat{K}_\gamma(n) := \frac{(K(n))^{\gamma}}{c(\gamma)},$$

so that $\sum_{n \in \mathbb{N}} \hat{K}_\gamma(n) = 1$ by the definition of $c(\gamma)$. It is important to realize that $\hat{K}_\gamma(\cdot)$ is still of the form (2.1), just with $\alpha$ replaced by $(1+\alpha)\gamma - 1 \geq 0$ and $L(\cdot)$ by $L(\cdot)^{\gamma}/c(\gamma)$ (which is still slowly varying at infinity).

Finally we need the following definitions: for $a, b, \nu \in \mathbb{R}$, $k \in \mathbb{N}$ and $0 < \gamma \leq 1$,

$$(2.10) \qquad f_\gamma(a, b; k) := \mathbb{E}\left[\left(\frac{1 + e^{a \sum_{i=1}^{k} \widetilde{\omega}_i + bk}}{2}\right)^{\gamma}\right]$$

and

$$
\begin{aligned}
G_\gamma(\nu, a, b) \\
:= \lim_{N\to\infty} \frac{1}{N} \log \hat{\mathbf{E}}_\gamma\left[ e^{\nu I_N(\tau)} \prod_{j=1}^{I_N(\tau)} f_\gamma(a, b; \tau_j - \tau_{j-1}) \mathbf{1}_{\{N \in \tau\}} \right],
\end{aligned}
\tag{2.11}
$$

where $\hat{\mathbf{P}}_\gamma$ is the law of the (recurrent) renewal with inter-arrival law $\hat{K}_\gamma(\cdot)$. The limit exists by superadditivity and is nonnegative.

Our main result is the following:

THEOREM 2.1. *Let $1/(1+\alpha) \leq \gamma \leq 1$ be such that $c(\gamma) < \infty$. Then*

$$F(\beta, h, \lambda, \widetilde{h}) \leq \frac{1}{\gamma} G_\gamma(\log c(\gamma) + h\gamma + \log M(\beta\gamma), -2\lambda, -2\lambda\widetilde{h}). \tag{2.12}$$

Some interesting consequences of this result are worked out in Section 3.

PROOF OF THEOREM 2.1. We have the elementary:

LEMMA 2.2 ([22], Chapter 2.1). *Let $0 < \gamma < 1$. If $n \in \mathbb{N}$ and $a_1 > 0, \ldots, a_n > 0$, then*

$$(a_1 + \cdots + a_n)^\gamma < a_1^\gamma + \cdots + a_n^\gamma. \tag{2.13}$$

We need also the identity

$$
\begin{aligned}
Z_{N,\omega,\widetilde{\omega}} \\
= \sum_{\ell=1}^{N} \sum_{0=i_0<i_1<\cdots<i_\ell=N} \left(\prod_{j=1}^{\ell} K(i_j - i_{j-1})\right) \exp\left(\beta \sum_{j=1}^{\ell} \omega_{i_j} + h\ell\right) \\
\times \prod_{j=1}^{\ell}\left(\frac{1 + \exp[-2\lambda \sum_{n=i_{j-1}+1}^{i_j}(\widetilde{\omega}_n + \widetilde{h})]}{2}\right),
\end{aligned}
\tag{2.14}
$$

which is just a way of rewriting the expectation in (2.2) as a sum over all possible configurations of $\tau \cap \{1, \ldots, N\}$. As a consequence of Lemma 2.2,

$$
\begin{aligned}
\mathbb{E}(Z_{N,\omega,\widetilde{\omega}})^\gamma \leq \mathbb{E}\Bigg[\sum_{\ell=1}^{N} \sum_{0=i_0<i_1<\cdots<i_\ell=N} \left(\prod_{j=1}^{\ell} K(i_j - i_{j-1})^\gamma\right) \\
\times \exp\left(\beta\gamma \sum_{j=1}^{\ell} \omega_{i_j} + h\gamma\ell\right) \\
\times \prod_{j=1}^{\ell}\left(\frac{1 + \exp[-2\lambda \sum_{n=i_{j-1}+1}^{i_j}(\widetilde{\omega}_n + \widetilde{h})]}{2}\right)^\gamma\Bigg].
\end{aligned}
\tag{2.15}
$$

Performing the disorder average, one sees that the right-hand side of (2.15) equals

$$
\sum_{\ell=1}^{N} \sum_{0=i_0<i_1<\cdots<i_\ell=N} \left(\prod_{j=1}^{\ell} \hat{K}_\gamma(i_j-i_{j-1})\right) e^{\ell[\log M(\gamma\beta)+h\gamma+\log c(\gamma)]}
$$
$$
\times \prod_{j=1}^{\ell} f_\gamma(-2\lambda, -2\lambda\widetilde{h}; i_j - i_{j-1}).
$$
(2.16)

Therefore,

$$
\limsup_{N\to\infty} \frac{1}{N} \log \mathbb{E}(Z_{N,\omega,\widetilde{\omega}})^\gamma
$$
$$
\leq G_\gamma(\log c(\gamma) + h\gamma + \log M(\beta\gamma), -2\lambda, -2\lambda\widetilde{h}).
$$
(2.17)

On the other hand for $\gamma > 0$ we have via Jensen's inequality:

(2.18) $\quad F(\beta, h, \lambda, \widetilde{h}) = \lim_{N\to\infty} \frac{1}{N} \mathbb{E} \log Z_{N,\omega,\widetilde{\omega}} \leq \limsup_{N\to\infty} \frac{1}{N\gamma} \log \mathbb{E}(Z_{N,\omega,\widetilde{\omega}})^\gamma$

which concludes the proof of Theorem 2.1. $\square$

## 3. Applications to disordered pinning models and copolymers.

3.1. *Random pinning model.* In this section we assume that $\lambda = 0$, $\mathbb{E}\omega_1 = 0$ and we call for simplicity the free energy $F(\beta, h)$. Since $\lambda = 0$ the random variables $\widetilde{\omega}_n$ do not play any role, and for simplicity we write $Z_{N,\omega}$ instead of $Z_{N,\omega,\widetilde{\omega}}$ for the partition function.

The model thus obtained is then the one considered for instance in [2, 12, 14, 18, 29]. According to the law **P** and, especially, to the value of $\alpha$, the model is interpreted in the physics literature as a pinning, or wetting, or Poland–Scheraga model in different spatial dimensions.

We observe that $F(\beta, \cdot)$ is nondecreasing and we denote as usual by $h_c(\beta)$ the (quenched) critical point of the pinning model:

(3.1) $\qquad h_c(\beta) := \inf\{h \in \mathbb{R} : F(\beta, h) > 0\},$

while the function $\beta \mapsto h_c(\beta)$ will be referred to as the *critical curve*. When $\beta = 0$ (homogeneous pinning model) it is a standard fact that $h_c(0) = -\log \mathbf{P}(\tau_1 < +\infty)$ (cf., for instance, [17], Chapter 2): $F(0, h)$ is positive for $h > h_c(0)$ and zero otherwise (for the detailed behavior of $F(0, h)$ for $h \searrow h_c(0)$, cf. [17], Theorem 2.1). It is also known (see [3] and [17], Chapter 5.2) that $h_c(\beta) < h_c(0)$ for every $\beta > 0$.

The annealed bound (2.6) applied to this case shows that

$$h_c(\beta) \geq h_c^{ann}(\beta) := -\log M(\beta) - \log \mathbf{P}(\tau_1 < +\infty) \tag{3.2}$$
$$= -\log M(\beta) - \log c(1).$$

On the other hand, since $f_\gamma(0,0;k) = 1$, Theorem 2.1 implies immediately

THEOREM 3.1. *For every $\beta > 0$ one has*

$$h_c(\beta) \geq \hat{h}_c(\beta) := \sup_{1/(1+\alpha) \leq \gamma \leq 1} -\frac{1}{\gamma} \log[M(\gamma\beta)c(\gamma)]. \tag{3.3}$$

Of course, $\hat{h}_c(\cdot)$ depends on $K(\cdot)$ because $c(\gamma)$ does.

The important point is that the bound provided by Theorem 3.1 is in various cases strictly better than the annealed one. For instance:

COROLLARY 3.2. *Assume that $\log M(\beta) \overset{\beta \to +\infty}{\sim} a\beta^\rho$ for some $\rho > 1$ and $a > 0$ [where $A(x) \overset{x \to \infty}{\sim} B(x)$ means $\lim_{x\to\infty} A(x)/B(x) = 1$].*

*Then, there exists $\beta_0 < \infty$ such that for every $\beta > \beta_0$*

$$h_c(\beta) > h_c^{ann}(\beta). \tag{3.4}$$

This applies for instance to the centered Gaussian case $\omega_1 \sim \mathcal{N}(0,1)$ where $\log M(\beta) = \beta^2/2$.

PROOF OF COROLLARY 3.2. Choose a value of $\gamma \in (1/(1+\alpha), 1)$. One has for $\beta \to \infty$:

$$-\frac{1}{\gamma} \log[M(\gamma\beta)c(\gamma)] \sim -a\gamma^{\rho-1}\beta^\rho \tag{3.5}$$

while $h_c^{ann}(\beta) \sim -a\beta^\rho$, and the statement is obvious from the conditions $\rho > 1$ and $\gamma < 1$. □

One can easily extract from Theorem 3.1 a *sufficient condition* which guarantees that $h_c(\beta) > h_c^{ann}(\beta)$ for a given $\beta > 0$, that is,

$$\partial_\gamma \left(-\frac{1}{\gamma} \log[M(\gamma\beta)c(\gamma)]\right)\bigg|_{\gamma=1} < 0. \tag{3.6}$$

For instance, in the case of Gaussian disorder and recurrent renewal [i.e., $c(1) = 1$], this condition can be re-expressed in the simple form:

$$\frac{\beta^2}{2} > -\sum_{n=1}^\infty K(n) \log K(n). \tag{3.7}$$

In the recent work [15] it is claimed, on the basis of numerical simulations, that in the case of the $(1+1)$-dimensional wetting model [which corresponds to $\alpha = 1/2$, $L(\cdot) \sim const$ and $c(1) = 1/2$] with $\omega_n$ taking only two values with equal probability, quenched and annealed critical points coincide even for strong disorder. Theorem 3.1 does not disprove this assertion because for symmetric two-valued $\omega_n$'s it turns out that equation (3.3) does not improve the annealed bound, but in our opinion it makes the scenario suggested by [15] rather unlikely.

REMARK 3.3. Corollary 3.2 is particularly interesting when $0 < \alpha < 1/2$ and disorder is Gaussian. Indeed, together with the results of [2] or [29], it implies that there exist $0 < \beta_1 \leq \beta_2 < \infty$ such that $h_c(\beta)$ coincides with the annealed critical point for $\beta \leq \beta_1$ and differs from it for $\beta > \beta_2$. In this sense, one can say that a transition from a weak disorder regime to a strong disorder regime occurs.

3.1.1. *Strong-disorder asymptotics of the critical point.* It is clear that, under the assumption of Corollary 3.2 on $M(\beta)$, for $\beta$ very large the supremum in (3.3) is realized by some $\gamma$ very close to $1/(1+\alpha)$. This observation, together with a generalization of ideas from [5], allows to identify the strong-disorder asymptotic behavior of the critical point. For instance, in the Gaussian case one has:

$$(3.8) \qquad h_c(\beta) \stackrel{\beta \to \infty}{\sim} -\frac{\beta^2}{2(1+\alpha)}$$

and it is easy to obtain an analogous statement in the case $\log M(\beta) \stackrel{\beta \to \infty}{\sim} a\beta^\rho$.

PROOF OF EQUATION (3.8). It is immediate to see from (3.3) that

$$(3.9) \qquad \liminf_{\beta \to \infty} \frac{h_c(\beta)}{\beta^2} \geq -\frac{1}{2(1+\alpha)}.$$

As for the opposite bound, it is based on a straightforward generalization of the *rare stretch strategy* of [5] (for this reason, we just sketch the main steps of the proof). Let $\ell \in \mathbb{N}$, assume that $N$ is an integer multiple of $\ell$ and divide $\{1, 2, \ldots, N\}$ into blocks $I_k := \{(k-1)\ell + 1, (k-1)\ell + 2, \ldots, k\ell\}$, with $k = 1, 2, \ldots, (N/\ell)$. Given $q > 0$ and the disorder realization $\omega$, let

$$(3.10) \qquad \mathcal{I}_\omega := \left\{ 1 \leq j \leq (N/\ell) : \sum_{n \in I_j} \omega_n \geq \ell q \right\} \cup \{N/\ell\}.$$

One obtains a lower bound on the partition function as follows:

$$Z_{N,\omega} \geq \mathbf{E}[e^{\sum_{n=1}^{N}(\beta\omega_n+h)\mathbf{1}_{\{n\in\tau\}}};$$
(3.11)
$$I_k \subset \tau \ \forall k \in \mathcal{I}_\omega; \tau \cap I_k = \varnothing \ \forall k \notin \mathcal{I}_\omega, k < (N/\ell)].$$

In other words, we have constrained $\tau$ to visit each point in the last block and in each of the blocks where the empirical average of the $\omega_n$'s is larger than $q$, and to skip all the others. We can now take the logarithm, divide by $N$ and let $N \to \infty$ at $\ell$ fixed in (3.11) to get a lower bound on the free energy. We do not detail this step, since an essentially identical computation appears in the proofs of [5], Proposition 3.1, [18], Theorem 3.1 and [17], Theorem 6.5. The net result is that for every $\varepsilon > 0$

$$(3.12) \quad F(\beta,h) \geq p(\ell)\left[\beta q + h + \log K(1) - (1+\alpha+\varepsilon)\frac{q^2}{2} + o_\ell(1)\right],$$

where

$$(3.13) \quad p(\ell) := \mathbb{P}\left(\sum_{n=1}^{\ell}\omega_n \geq \ell q\right)$$

and $o_\ell(1)$ is a quantity which vanishes for $\ell \to \infty$. The term $\log K(1)$ is due to $\mathbf{P}(I_1 \subset \tau) = K(1)^\ell$.

From (3.12) and the definition of the critical point one deduces that

$$(3.14) \quad h_c(\beta) \leq -\beta q + (1+\alpha+\varepsilon)\frac{q^2}{2} - \log K(1).$$

Optimizing over $q$ and using the arbitrariness of $\varepsilon > 0$,

$$(3.15) \quad h_c(\beta) \leq -\frac{\beta^2}{2(1+\alpha)} - \log K(1)$$

which, together with (3.9), proves equation (3.8). $\square$

3.1.2. *About the size of $Z_{N,\omega}$ in the delocalized phase.* For the homogeneous model, it is known [17], Theorem 2.2, that if $h < h_c(0)$ then

$$(3.16) \quad Z_N = \mathbf{E}[e^{hI_N(\tau)}\mathbf{1}_{\{N\in\tau\}}] \overset{N\to\infty}{\sim} CK(N)$$

where $C > 0$ depends on $h$ and on $\mathbf{P}(\tau_1 < \infty)$. It is natural to ask whether a similar statement holds for the disordered model inside the delocalized phase. In this respect, the ideas developed in this work allow to go much beyond the statement of Theorem 3.1 that the infinite-volume free energy vanishes for $h \leq \hat{h}_c(\beta)$, and prove the following: if $h < \hat{h}_c(\beta)$ there exists $1/(1+\alpha) < \gamma \leq 1$ and a constant $C := C(h, K(\cdot)) < \infty$ such that

$$(3.17) \quad \mathbb{P}(Z_{N,\omega} \geq uK(N)) \leq Cu^{-\gamma}$$

for every $u > 0$. The upper bound (3.17) on the partition function should be read together with the lower bound (2.5): both are of order $K(N)$, just like the estimate (3.16) which holds for the pure model.

PROOF OF EQUATION (3.17). Since $h < \hat{h}_c(\beta)$, there exists $1/(1+\alpha) < \gamma \leq 1$ such that

$$\Delta := -\log[M(\beta\gamma)c(\gamma)] - h\gamma > 0.$$

Then, it follows from the proof of Theorem 2.1 [cf. in particular equations (2.15) and (2.16) taken for $\lambda = 0$] and from (3.16) that

$$(3.18) \quad \mathbb{E}(Z_{N,\omega})^\gamma \leq \hat{\mathbf{E}}_\gamma[e^{-I_N(\tau)\Delta}\mathbf{1}_{\{N\in\tau\}}] \leq C\hat{K}_\gamma(N) = \frac{C}{c(\gamma)}K(N)^\gamma,$$

for some $C := C(h, K(\cdot))$. Equation (3.17) is then an immediate consequence of Markov's inequality. □

Let us remark also that one can extract from the proof of (3.17) the following almost sure statement: if $h < -(1/\gamma)\log[M(\beta\gamma)c(\gamma)]$ for some $1/(1+\alpha) < \gamma \leq 1$, then

$$(3.19) \qquad \frac{Z_{N,\omega}}{K(N)^\mu} \stackrel{N\to\infty}{\longrightarrow} 0,$$

$\mathbb{P}(d\omega)$-almost surely, for every $\mu < 1 - 1/(\gamma(1+\alpha))$.

3.2. *"Reduced" wetting model.* This model, introduced in [5], Section 4 as a toy version of the copolymer model discussed in the next section, is obtained from (2.2) putting $\lambda = h = 0$ and assuming that $\mathbb{P}(\omega_1 = 1) = p = 1 - \mathbb{P}(\omega_1 = 0)$, for some $p \in (0,1)$. Moreover, one assumes here that the renewal is transient, $c(1) < 1$, otherwise the model is uninteresting in that no phase transition occurs. As we will see, the actual value of $c(1)$ is irrelevant in the strong-disorder regime we are going to consider, as long as it is strictly smaller than 1. Changing conventions with respect to the previous section, we denote the free energy as $F(\beta, p)$ in this case.

It is known that for every $\beta > 0$ there exists $p_c(\beta) \in (0,1)$ such that the free energy $F(\beta, p)$ is zero for $p \leq p_c(\beta)$ and positive for $p > p_c(\beta)$. The main question is to compute (if it exists) the limit

$$(3.20) \qquad m_c := -\lim_{\beta\to\infty}\frac{1}{\beta}\log p_c(\beta).$$

It was proven in [5] that for $\alpha = 1/2$

$$(3.21) \qquad 2/3 \leq -\limsup_{\beta\to\infty}\frac{1}{\beta}\log p_c(\beta) \leq -\liminf_{\beta\to\infty}\frac{1}{\beta}\log p_c(\beta) \leq 1.$$

This result can be easily generalized to any $\alpha > 0$ and in this case the lower bound in (3.21) is replaced by $1/(1+\alpha)$.

Here we prove the following:

THEOREM 3.4. *The limit in* (3.20) *exists and equals* $1/(1+\alpha)$.

Recently a proof of Theorem 3.4, based on entirely different ideas, was given in [6]. While our approach is much simpler, the method of Bolthausen, Caravenna and de Tilière is more natural from a renormalization group point of view.

PROOF OF THEOREM 3.4. We start by remarking that if we apply Theorem 2.1 to the reduced model we obtain immediately

$$(3.22) \qquad -\liminf_{\beta \to \infty} \frac{1}{\beta} \log p_c(\beta) < 1,$$

but not the sharper statement of Theorem 3.4. To go beyond (3.22), somehow we have to use the information that, if $p < p_c(\beta)$, the sites where $\omega_n = 1$ are very sparse for $\beta$ large (their density $p$ is indeed exponentially small in $\beta$), and that between two such sites $\tau$ has typically just a finite number of points, since it is a transient renewal under the law $\mathbf{P}$.

Let $1/(1+\alpha) < \gamma < 1$ and $C > 1$ such that

$$(3.23) \qquad \frac{c(\gamma)}{C^\gamma} < 1.$$

The renewal $\tau$ with inter-arrival law $K(\cdot)$ being transient, that is, $c(1) < 1$, it is known (cf., for instance, [17], Theorem A.4) that

$$(3.24) \qquad \mathbf{P}(n \in \tau) \overset{n \to \infty}{\sim} \frac{1}{(1-c(1))^2} K(n),$$

and therefore there exists $C' := C'(\gamma, C) < \infty$ such that

$$(3.25) \qquad \mathbf{P}(n \in \tau) \leq C' \frac{K(n)}{C}$$

for every $n \in \mathbb{N}$.

We let now $\mathcal{Y}_\omega := \{1 \leq n \leq N : \omega_n = 1\} \cup \{N\}$ and we decompose the partition function according to the configuration of $A(\tau) := \tau \cap \mathcal{Y}_\omega$. If $|A(\tau)| = \ell (\geq 1)$, we write $A(\tau) = \{a_1, a_2, \ldots, a_\ell\}$ with $a_i < a_{i+1}$ and, by convention, we set $a_0 := 0$. Then,

$$(3.26) \qquad Z_{N,\omega} \leq \sum_{\ell=1}^{N} \sum_{\substack{A \subseteq \mathcal{Y}_\omega: \\ |A|=\ell, a_\ell=N}} e^{\beta \ell} \prod_{j=1}^{\ell} \mathbf{P}(a_j - a_{j-1} \in \tau)$$

$$(3.27) \quad \leq \sum_{\ell=1}^{N} \sum_{\substack{A \subseteq \mathcal{Y}_\omega: \\ |A|=\ell, a_\ell=N}} e^{(\beta+\log C')\ell} \prod_{j=1}^{\ell} \frac{K(a_j - a_{j-1})}{C}$$

$$\leq e^{\beta+\log C'} \sum_{\ell=1}^{N} \sum_{0=a_0<a_1<\cdots<a_\ell=N} e^{(\beta+\log C')\sum_{j=1}^{\ell} \omega_{a_j}}$$

$$(3.28) \quad \times \prod_{j=1}^{\ell} \frac{K(a_j - a_{j-1})}{C},$$

where in the first inequality we used the fact that

$$(3.29) \quad \begin{aligned} \mathbf{P}(a_i \in \tau, \tau \cap \mathcal{Y}_\omega \cap \{a_{i-1}+1, \ldots, a_i - 1\} = \varnothing | a_{i-1} \in \tau) \\ \leq \mathbf{P}(a_i - a_{i-1} \in \tau), \end{aligned}$$

in the second one we used (3.25) and the third one is obvious since we have just added extra positive terms to the sum. On the other hand the right-hand side of (3.28), apart from the global factor $\exp(\beta + \log C')$ which is anyway irrelevant for the computation of the infinite-volume limit of the free energy, is just the partition function of the model where $\beta$ is replaced by $(\beta + \log C')$ and $K(n)$ by $K(n)/C$ for every $n \in \mathbb{N}$ [cf. equation (2.14) taken for $\lambda = 0$]. Following step by step the proof of Theorem 2.1, one finds therefore that

$$\lim_{N \to \infty} \frac{1}{N} \mathbb{E} \log Z_{N,\omega}$$

$$(3.30) \quad \leq \limsup_{N \to \infty} \frac{1}{N\gamma} \log \mathbb{E}(Z_{N,\omega})^\gamma$$

$$\leq \lim_{N \to \infty} \frac{1}{N\gamma} \log \hat{\mathbf{E}}_\gamma[e^{I_N(\tau) \log[(c(\gamma)/C^\gamma)(pe^{\gamma(\beta+\log C')}+(1-p))]} \mathbf{1}_{\{N \in \tau\}}],$$

which vanishes whenever

$$(3.31) \quad \log[pe^{\gamma(\beta+\log C')} + (1-p)] + \log\left(\frac{c(\gamma)}{C^\gamma}\right) \leq 0.$$

Thanks to the choice (3.23), it follows immediately that

$$(3.32) \quad -\liminf_{\beta \to \infty} \frac{1}{\beta} \log p_c(\beta) \leq \gamma.$$

By the arbitrariness of $(1 \geq)\gamma > 1/(1+\alpha)$, and since we already know that

$$-\limsup_{\beta \to \infty} \frac{1}{\beta} \log p_c(\beta) \geq 1/(1+\alpha),$$

we obtain the statement of the theorem. □

3.3. *Copolymer model.* In this section we set $\beta = h = 0$ and with yet another abuse of notation we call the free energy $F(\lambda, \widetilde{h})$. This is just (a generalization of) the copolymer model considered for instance in [4, 5, 7].

REMARK 3.5. The results which follow can be easily generalized to the case $h \neq 0$. This is particularly evident for $h < 0$. Indeed, it is well known [and immediate to check from (2.14)] that the model with $h < 0$ is equivalent to the one with $h = 0$ and $K(n)$ replaced by $K(n) \exp(-|h|)$ for every $n \in \mathbb{N}$ [so that $c(\gamma)$ becomes $c(\gamma) \exp(-\gamma |h|)$].

We observe that in view of $\lambda \geq 0$ the free energy is nonincreasing with respect to $\widetilde{h}$ and we denote by $\widetilde{h}_c(\lambda)$ the quenched critical point:

$$(3.33) \qquad \widetilde{h}_c(\lambda) := \sup\{\widetilde{h} \geq 0 : F(\lambda, \widetilde{h}) > 0\},$$

while $\lambda \mapsto \widetilde{h}_c(\lambda)$ is the critical line.

It is convenient to define for $0 < \gamma \leq 1$

$$(3.34) \qquad \widetilde{h}_c^{(\gamma)}(\lambda) := \frac{1}{2\lambda\gamma} \log \widetilde{M}(-2\lambda\gamma).$$

Note that $\widetilde{h}_c^{(\cdot)}(\lambda)$ is increasing by Jensen's inequality and that

$$(3.35) \qquad \lim_{\lambda \searrow 0} \frac{\widetilde{h}_c^{(\gamma)}(\lambda)}{\lambda} = \gamma$$

thanks to $\mathbb{E}\widetilde{\omega}_1^2 = 1$, $\mathbb{E}\widetilde{\omega}_1 = 0$. For $\gamma = 1$, this is just the annealed critical line: $F^{ann}(0, 0, \lambda, \widetilde{h}) = 0$ if and only if $\widetilde{h} \geq \widetilde{h}_c^{(1)}(\lambda)$, as it is immediate to realize from (2.6). On the other hand, $\widetilde{h}_c^{(1/(1+\alpha))}(\cdot)$ is sometimes referred to as the "Monthus line." This line was proposed as the true critical line of the copolymer model in the theoretical physics literature, on the basis of renormalization group arguments [25].

In general one knows that for every $\lambda > 0$,

$$(3.36) \qquad \widetilde{h}_c^{(1/(1+\alpha))}(\lambda) \leq \widetilde{h}_c(\lambda) \leq \widetilde{h}_c^{(1)}(\lambda).$$

The lower bound in (3.36) was proven in [5] in the case where $K(\cdot) = K^{SRW}(\cdot)$, the law of the first return to zero of the one-dimensional simple random walk (cf. Section 2), but the bound can be proven to hold for every $K(\cdot)$ of the form (2.1), see [17], Chapter 6. Numerical simulations (supported by solid probabilistic estimates) performed *in the particular case* $K(\cdot) = K^{SRW}(\cdot)$ strongly indicate that neither of the two bounds in (3.36) is optimal [9].

For $K(\cdot) = K^{SRW}(\cdot)$, it is also known that the limit slope $\lim_{\lambda \to 0} \widetilde{h}_c(\lambda)/\lambda$ exists [7] [and is therefore bounded between $2/3$ and $1$ in view of (3.36) and $\alpha = 1/2$], that it coincides with the slope of the critical curve of a continuous copolymer model with Brownian disorder [7] and that it is to a large extent independent of $\mathbb{P}$ [19].

Here we will show that (say, in the Gaussian case) the annealed upper bound $\widetilde{h}_c(\lambda) \leq \widetilde{h}_c^{(1)}(\lambda)$ can be improved for large $\lambda$ whatever $K(\cdot)$ is [within the class (2.1)]. If the renewal is transient, $\sum_{n \in \mathbb{N}} K(n) < 1$, then with no assumptions on the disorder distribution annealing can be improved for every $\lambda > 0$ and, in particular, the slope at the origin turns out to be strictly smaller than 1. Finally, if $K(\cdot)$ satisfies condition (3.43) below, the lower bound in (3.36) is the optimal one for every $\lambda$. This condition *is not satisfied* for $K(\cdot) = K^{SRW}(\cdot)$, in agreement with the simulations mentioned above. By the way, our results strongly indicate that the critical curve depends in general not only on the tail behavior of $K(\cdot)$ (say, on the exponent $\alpha$), as one might be tempted to guess on the basis of belief in universality, but also on the details of the slowly-varying function $L(\cdot)$. Note that this is not the case for the annealed curve, which depends only on the disorder distribution $\mathbb{P}$, nor for the Monthus line which depends only on $\mathbb{P}$ and $\alpha$.

Observe first of all that, thanks to Lemma 2.2,

$$(3.37) \qquad f_\gamma(a,b;k) \leq 2^{1-\gamma} \frac{1 + e^{k(b\gamma + \log \widetilde{M}(\gamma a))}}{2}.$$

Therefore,

$$(3.38) \quad \begin{aligned} G_\gamma(\nu, a, b) \\ \leq \lim_{N \to \infty} \frac{1}{N} \log \hat{\mathbf{E}}_\gamma \Bigg[ & e^{I_N(\tau)(\nu + (1-\gamma) \log 2)} \\ & \times \prod_{j=1}^{I_N(\tau)} \left( \frac{1 + e^{(b\gamma + \log \widetilde{M}(\gamma a))(\tau_j - \tau_{j-1})}}{2} \right) \mathbf{1}_{\{N \in \tau\}} \Bigg]. \end{aligned}$$

Our main result for the copolymer model is the following:

THEOREM 3.6. *Define*

$$(3.39) \qquad \gamma_c := \inf\{0 < \gamma \leq 1 : \log c(\gamma) < \gamma \log 2\} \in [1/(1+\alpha), 1).$$

*Then, one has:*

1. *For every $\gamma > \gamma_c$ there exists $C(\gamma) < \infty$ such that for every $\lambda > 0$*

$$(3.40) \qquad \widetilde{h}_c(\lambda) \leq \widetilde{h}_c^{(\gamma)}(\lambda) + \frac{C(\gamma)}{\lambda}.$$

2. *If $\gamma_c = 1/(1+\alpha)$ and*

(3.41) $$0 < \log c(1/(1+\alpha)) + \frac{\alpha}{1+\alpha} \log 2 < \log 2$$

*there exists $C < \infty$ such that for every $\lambda > 0$*

(3.42) $$\widetilde{h}_c(\lambda) \leq \widetilde{h}_c^{(1/(1+\alpha))}(\lambda) + \frac{C}{\lambda}.$$

3. *If $\gamma_c = 1/(1+\alpha)$ and*

(3.43) $$\log c(1/(1+\alpha)) + \frac{\alpha}{1+\alpha} \log 2 \leq 0,$$

*then for every $\lambda > 0$*

(3.44) $$\widetilde{h}_c(\lambda) = \widetilde{h}_c^{(1/(1+\alpha))}(\lambda).$$

In view of Remark 3.5 above, condition (3.43) is realized for instance if $c(1/(1+\alpha)) < \infty$ and we add a homogeneous depinning term proportional to $-|h|$, with $|h|$ sufficiently large. In any case, we emphasize that a necessary (but not sufficient) condition for (3.43) to hold is that $\tau$ be transient under **P**.

REMARK 3.7. Note that in the case $L(\cdot) \sim const$, for example, if $K(\cdot) = K^{SRW}(\cdot)$, one has $\gamma_c > 1/(1+\alpha)$ since $c(\gamma) \to \infty$ for $\gamma \searrow 1/(1+\alpha)$. In particular, from the explicit expression of $K^{SRW}(\cdot)$ one finds numerically $\gamma_c \simeq 0.83$ in the simple random walk case. It is interesting to note that the numerical results of [9], Section 5, are compatible with $\widetilde{h}_c(\lambda) = \widetilde{h}_c^{(m)}(\lambda)$ with $m$ somewhere between 0.8 and 0.84. Understanding whether this is more than just a coincidence requires further numerical simulations of the type [9], performed also for other choices of $K(\cdot)$.

We have also the analogue of Corollary 3.2:

COROLLARY 3.8. *Assume that $\log \widetilde{M}(-\lambda) \sim a\lambda^\rho$ for $\lambda \to +\infty$, for some $a > 0$ and $\rho > 1$. Then, there exists $\lambda_0 < \infty$ such that for $\lambda > \lambda_0$*

(3.45) $$\widetilde{h}_c(\lambda) < \widetilde{h}_c^{(1)}(\lambda).$$

This is obvious from (3.40) since $\gamma_c < 1$.

Finally, as we mentioned, in the transient case (and without assumptions on $\mathbb{P}$) the improvement on annealing can be pushed down to $\lambda = 0$:

COROLLARY 3.9. *Assume that $c(1) = \mathbf{P}(\tau_1 < +\infty) < 1$. Then, for every $\lambda > 0$*

(3.46) $$\widetilde{h}_c(\lambda) \leq \widetilde{h}_c^{(\bar\gamma)}(\lambda)$$

*where*

(3.47) $$\bar\gamma := \inf\{\gamma \leq 1 : \log c(\gamma) + (1-\gamma) \log 2 \leq 0\} \in [1/(1+\alpha), 1).$$

As a consequence of equation (3.35) we have that if $\mathbf{P}(\tau_1 < +\infty) < 1$ (or if $h < 0$, see Remark 3.5) the slope of the critical curve at the origin (if it exists) is strictly smaller than 1.

PROOF OF THEOREM 3.6. One has from (3.38)

$$G_\gamma(\log c(\gamma), -2\lambda, -2\lambda[\widetilde{h}_c^{(\gamma)}(\lambda) + \varepsilon])$$

$$(3.48) \qquad \leq \lim_{N\to\infty} \frac{1}{N} \log \hat{\mathbf{E}}_\gamma \left[ e^{I_N(\tau)(\log c(\gamma) + (1-\gamma)\log 2)} \right.$$

$$\left. \times \prod_{j=1}^{I_N(\tau)} \left( \frac{1 + e^{-2\lambda\gamma\varepsilon(\tau_j - \tau_{j-1})}}{2} \right) \mathbf{1}_{\{N \in \tau\}} \right].$$

Let us consider first case (1). To this purpose, let $\gamma > \gamma_c$ and choose $\varepsilon = C/(\gamma\lambda)$ in (3.48). It is clear that, for every $\delta > 0$, it is possible to choose $C = C(\delta)$ sufficiently large so that

$$(3.49) \qquad \frac{1 + e^{-2Ck}}{2} \leq e^{-\log 2 + \delta}$$

for every $k \in \mathbb{N}$. Taking for instance

$$(3.50) \qquad \delta = \delta(\gamma) = \frac{\gamma \log 2 - \log c(\gamma)}{2} > 0$$

(so that $C$ depends on $\gamma$) one finds therefore

$$(3.51) \qquad G_\gamma\left(\log c(\gamma), -2\lambda, -2\lambda\left[\widetilde{h}_c^{(\gamma)}(\lambda) + \frac{C(\gamma)}{\gamma\lambda}\right]\right)$$

$$\leq \lim_{N\to\infty} \frac{1}{N} \log \hat{\mathbf{E}}_\gamma[e^{-I_N(\tau)\delta(\gamma)} \mathbf{1}_{\{N \in \tau\}}] = 0.$$

[The positivity of $\delta(\gamma)$ follows from $\gamma > \gamma_c$ and from the definition of $\gamma_c$.] By Theorem 2.1 one has therefore (3.40). Equation (3.42) is proven analogously.

Finally, if condition (3.43) is verified then choosing $\gamma = 1/(1+\alpha)$ and $\varepsilon = 0$ in equation (3.48) one has

$$(3.52) \qquad G_{1/(1+\alpha)}\left(\log c\left(\frac{1}{1+\alpha}\right), -2\lambda, -2\lambda\widetilde{h}_c^{(1/(1+\alpha))}(\lambda)\right)$$

$$\leq \lim_{N\to\infty} \frac{1}{N} \log \hat{\mathbf{E}}_{1/(1+\alpha)}[e^{I_N(\tau)(\log c(1/(1+\alpha)) + \alpha/(1+\alpha)\log 2)} \mathbf{1}_{\{N \in \tau\}}]$$

and the right-hand side is zero since it is the free energy of a homogeneous pinning model with nonpositive pinning strength, see assumption (3.43). This shows that $\widetilde{h}_c(\lambda) \leq \widetilde{h}_c^{(1/(1+\alpha))}(\lambda)$, and the opposite bound is already in (3.36). □

PROOF OF COROLLARY 3.9. Just take equation (3.48) for $\gamma = \bar{\gamma}$ and $\varepsilon = 0$. □

**Acknowledgments.** I wish to thank Erwin Bolthausen, Bernard Derrida and Giambattista Giacomin for interesting discussions.

ÉCOLE NORMALE SUPÉRIEURE DE LYON
LABORATOIRE DE PHYSIQUE AND CNRS
UMR 5672
46 ALLÉE D'ITALIE
69364 LYON CEDEX 07
FRANCE
E-MAIL: fltonine@ens-lyon.fr